\documentclass[english,french,12pt,twoside]{article}
\usepackage{latexsym,graphics,color,amssymb,babel}
\title{Weak amenability of Coxeter groups}
\author{Gero Fendler \\ 
Naturwissenschaftlich-Technische Fakult\"{a}t I\\%
Fachbereich 6.1 Mathematik \\
Universit\"at des Saarlandes \\ 
}
\newcounter{prop}
\newcounter{def}

\newcommand{\abs}[1][.]{|{#1}|}
\newcommand{\norm}[2][]{\| \, {#2} \,\|_{#1}}
\newcommand{\scalprod}[3][\alpha]{\int_{0}^{\infty}\int_{0}^{\infty}{{#1}^{s \wedge t} {#2}(s) \overline{{#3}(t)}} \, ds \, dt}
\newcommand{\scalproddis}[3][\alpha]{\sum_{k,l=0}^{\infty}{{#1}^{k \wedge l}{#2}_k \overline{{#3}_l}}}
\newcommand{\card}[1][\#]{\protect{#1}}
\newcommand{\e}{\mathbf{e}}
\newcommand{\di}[2]{{\rm \bf d}(#1,#2)}
\newcommand{\m}[2]{{\rm m}(#1,#2)}
\newcommand{\la}[1]{{\rm \bf l}(#1)}

\newcommand{\cayley}[2]{\mathcal{C}(#1,#2)}
\newlength{\expheight}
\newcommand{\heightheta}{\mbox{\raisebox{0pt}[\expheight][0pt]{$\theta$}}}
\newtheorem{thm}{Theorem}
\newtheorem{cor}{Corollary}
\newtheorem{lem}{Lemma}

\newtheorem{rem}{Remark}

\newtheorem{prop}[prop]{Proposition}

\newenvironment{proof}[1][Proof:]{\it #1 \rm }{}
\begin{document}
\selectlanguage{english}
\maketitle
\vspace{1cm}
\begin{abstract}
{Let $(G,S)$ be a Coxeter group. We construct a continuation,
to the open unit disc, of the unitary representations associated
to the positive definite functions $g\mapsto r^{\la{g}}$. (Here
$0<r<1$, and ${\rm \bf l}$ denotes the length function with respect to the
generating set $S$.)
\par
The constructed representations are uniformly bounded and we 
prove that this implies the weak amenability of the group $G$.} 
\end{abstract}
\selectlanguage{french}
\begin{abstract}
Soit $(G,S)$ un groupe de Coxeter. Nous construisons une extension,
au disque unit\'{e} ouvert, de la s\'{e}rie des repr\'{e}sentations unitaires
associ\'{e}es aux fonctions definies positives  $g\mapsto r^{\la{g}}$.
(Ici, $0<r<1$, et ${\rm \bf l}$ d\'{e}signe  la fonction longueur
par rapport aux g\'{e}n\'{e}rateurs $S$.)
\par
Les repr\'{e}sentations ainsi construites sont uniform\'{e}ment born\'{e}es et
nous d\'{e}\-montrons que ceci implique que le groupe $G$
est faiblement moyenable.
\end{abstract}
\selectlanguage{english}
\vfill
\thanks{{\bf Acknowledgements :} Part of this work has been done while
the author visited the Mathematical Institute of the University of Wroc{\l}aw
in 1998.}
\section{Introduction}
Let $G$ be a locally compact group. On the space
$L^2(G)$ of, with respect to a left invariant Haar measure,
square integrable functions the left regular
representation of $G$ is defined by the unitary operators
$\lambda(g)\;,\;g\in G$, where those are given by
\[\lambda(g)f\,(h)=f(g^{-1}h)\quad h\in G,\,f\in L^2(G).\]
The reduced von~Neumann algebra $VN_r(G)\subset B(L^2(G))$
is the weak operator topology closure of the linear span of
$\{\lambda(g)\;:\;g\in G\,\}$.
\par
A complex valued function $\phi$ on $G$
is called positive definite
whenever for  finite subsets $g_1, \ldots ,g_n\in G$,
$c_1,\ldots,c_n \in \mathbb{C}$:
\[ \sum_{i,j=1}^n c_i\overline{c}_j\phi(g_{i}^{-1}g_{j}) \geq 0.\]
We denote $A(G)$ the Fourier algebra of $G$ as defined by 
Eymard~\cite{Eymard64} and recall that on one hand 
it is the linear span of the positive definite functions
$g\mapsto <\lambda(g)\; f, f>$, $f\in L^2(G)$, which actually is
an algebra of continuous functions on $G$, and on the other
hand is naturally identified with the predual of $VN_r(G)$.
\par
In terms of the Fourier algebra amenability of $G$ can be characterized
by the existence of a $A(G)$--norm bounded approximate
identity in $A(G)$. That is, there exist a constant $C>0$ and 
a net $(m_i)_{i\in I} \in A(G)$ such that
\begin{eqnarray*}
\lim_i \norm[A(G)]{m_i\varphi -\varphi} &=&0 \qquad \forall\; \varphi \in A(G)\\
\sup_i \norm[A(G)]{m_i} &=& C.
\end{eqnarray*}
\par
In their work on multipliers of the Fourier algebra of
some simple Lie groups and their discrete subgroups \cite{CanHaa_84}
de~Canni\`{e}re and Haagerup started to investigate more 
general approximate identities of the Fourier algebra
than norm bounded ones.
\par
They introduced the concept of completely bounded multipliers
in that paper:
\par
A function $m$ on the locally compact group $G$  
is a multiplier of the  Fourier algebra $A(G)$, if
for any $\varphi \in A(G)$
the pointwise product of functions $M_m(\varphi)=m\cdot\varphi\in A(G)$ again.
Now, the dual to this multiplication operator acts on 
the reduced von~Neumann algebra $VN_r(G)\subset B(L^2(G))$. From this
last inclusion one can let act the $n\times n $ matrices with entries 
from $VN_r(G)$ on $l^{2}_{n}\otimes_2L^2(G)$ in the canonical way
and norm $M_n(VN_r(G))$ accordingly, for each $n\in \mathbb{N}$.
\par
The action of a linear operator on $VN_r(G)$ can be extended
to a linear action on  $M_n(VN_r(G))$ simply by letting it act
on each matrix entry (of course we think of $M_{m}^{\ast}$).
The multiplier $m$ of $A(G)$ is called a completely bounded
multiplier of the Fourier algebra if those extensions
are bounded, uniformly  in $n\in \mathbb{N}$. Clearly, taking the smallest 
possible 
bound defines a norm $\norm[M_0A(G)]{m}$ on a subspace of all multipliers
of the Fourier algebra.
\par
A locally compact group is {\em weakly amenable}, 
if there exist a constant $C \; < \;\infty$
and a net $(m_i)_{i\in I} \in A(G)$ such that
\begin{eqnarray*}
\lim_i \norm[A(G)]{m_i\varphi -\varphi} &=&0 \qquad \forall\; \varphi \in A(G)\\
\sup_i \norm[M_0A(G)]{m_i} &=& C.
\end{eqnarray*}
The least such constant is called the Cowling Haagerup constant of $G$.
\par 
Since $\norm[M_0A(G)]{}\leq \norm[A(G)]{}$, an amenable group
is weakly amenable too, but on the other hand
there are non--amenable groups which are weakly amenable.
\par
It turned out \cite{BozeFend_84} that the space of completely bounded
multipliers of the Fourier algebra
had been discussed in harmonic analysis by Herz \cite{Herz_74} as
a generalisation of the Fourier Stieltjes algebra, the algebra of coefficients
of continuous unitary representations. 
In their above cited paper, de~Canni\`{e}re and Haagerup 
showed that the Fourier algebras of 
finite extensions of the Lie groups $SO_0(n,1), \; n\geq2$
indeed contain completely bounded approximate identities.
This result transfers to the discrete subgroups of those groups too.
A special example is the free group on two generators, $\mathbb{F}_2$.
\par
Subsequently, this result was extended on one hand to simple 
Lie groups of real rank one by Cowling and Haagerup \cite{CowHaa_89}
and, on the other hand, to amalgamated products by
Bo\.{z}ejko and Picardello~\cite{BozPic} and to groups acting on trees by 
Szwarc \cite{Szwarc_91} and Valette \cite{Valette_90}.
\par
The construction of such an approximate identity
relies on the existence of a path of uniformly bounded
representations which connect, in some way, the regular representation to the
trivial representation of the group in question. Moreover, and we see
no way around this difficulty, one has to extend this path,
in some way continuously, to some complex region, still preserving
the uniform boundedness of the representations.
We note that, for non-amenable groups, it is not possible to deal 
here only with
unitary representation. Whereas the original path of representations 
might consist of unitary ones. (Of course if the group has 
the Kazdhan property then this is not possible either.)
\par
In this paper we shall continue to consider some new examples
from the variety of discrete groups.
In fact we shall deal with Coxeter systems $(G,S)$,
which we, by abuse of language call Coxeter groups, the generating set
is always understood. 
\par
For Coxeter groups
it is known from the work of 
Bo\.{z}ejko, Januszkiewicz and Spatzier \cite{BozJanSpatz_88} that
the length function, with respect to $S$ is negative definite,
where a function $\phi : G \to \mathbb{C}$ is called negative definite, whenever
for finitely many $g_1,\ldots,g_n \in G$ and $c_1, \ldots,c_n \in \mathbb{C}$
with $\sum_{i=1}^n c_i =0$:
\[\sum_{i,j=1}^n c_i \overline{c}_j \phi(g_i^{-1} g_j) \leq 0. \]
Hence, by a theorem of Schoenberg (see e.\ g.~\cite{BergForst_75}), for $0<r< 1$, 
\[
g \mapsto r^{\la{g}} \qquad g\in G,
\]
is a positive definite function and the associated representations
form   
a suitable path of unitary representations.   
\par
Starting from this we consider the problem of extending this series of
representations to a complex parameter 
$z\in D=\{z\in \mathbb{C}\,:\, \abs[z]<1\}$.
\par
The author knows of a manuscript
of T. Januszkiewicz
in which it is proved that for all finitely generated
Coxeter groups the functions $g\mapsto z^{\la{g}}$, 
$z\in \mathbb{C} \it,\; 0\leq|z|<1$,
are coefficients of uniformly bounded representations \cite{Januszkiewicz_99}.
We learned a geometric Lemma from this, which is used in work of
Millson, see Lemma~2.1 of \cite{Millson_76}. Millson attributes it to Jaffe.
For the readers convenience we shall state it as Lemma~\ref{lem:disjoint} and prove it 
in a formulation convenient for us.
\par
The paper is organised as follows.
After this introduction
in section 2 we define some positive definite kernels
and discuss domination properties between them, 
which we shall apply to prove bounds on representations of the Coxeter group.
In section 3 we give an introduction to the standard geometrical 
representation of a Coxeter group and show that the geometry implies
the positive definiteness of certain kernels related to the length
function of the group. These results are then used
to prove uniform bounds on representations  constructed by
modifying the Gelfand-Naimark-Segal representation associated to
the positive definite functions 
$g\mapsto r^{\la{g}}$, $0<r<1$  in section 4.
Section 5 finally contains a proof of the weak amenability.
\section{Some analytical tools}
Denoting $\mathbb{R}_+$ the non--negative reals and
$s \wedge t= \min (s,t)\;s,t \in \mathbb{R}$ the minimum of two reals 
we define for $ \alpha > 1 $ a kernel
$ k_{\alpha} \, : \,  \mathbb{R}_+ \times \mathbb{R}_+ \rightarrow \,\mathbb{R}_+ $
by 
\[  k_{\alpha}(s,t) \, := \, \alpha ^{s \wedge t} \quad s,t \in \mathbb{R}_+ \]
and a corresponding sesquilinear form on the space $ \mathcal{F}(\mathbb{R}_+) $ of 
complex valued compactly supported locally Lebesgue integrable functions 
defined on  
$ \mathbb{R}_+ :$
\[ [ \, f \, , \, g \, ]_{\alpha} \, := \, \scalprod{f}{g}\, \quad 
f,g \in \mathcal{F}{(\mathbb{R}_+)}. \]
\par
It is then easy to see that the sesquilinear form is non-negative definite.
For, if 
$ s,t \, \geq \,0, $ 
then
\[ s \wedge t \, = \, \int_{0}^{\infty}\,{ \chi_{[0,s)} (v) \chi_{[0,t)} (v)} \, dv \, , \]
where
$ \chi_{[0,s)} $
denotes the indicator function of the interval 
$[0,s).$
Thus,
$ (s,t) \, \mapsto \, s \wedge t $
is a positive definite and hence
$ (s,t) \, \mapsto \, - s \wedge t $
is a negative definite kernel on 
$ \mathbb{R} _+.$
Since 
$ \log(\alpha) \, > \, 0 ,$ we obtain by Sch\"{o}nberg's 
theorem~\cite{BergForst_75} that
\[ k_{\alpha} \, : \, (s,t) \mapsto \alpha^{s \wedge t} \, = \,
e^{-\log(\alpha )(-s \wedge t)} \quad \quad s,t \in \mathbb{R}_{+} \]
is positive definite. Moreover, from its series development
\begin{eqnarray*}
  \label{eq:posdef}
  k_{\alpha}(s,t) & = & 1\, +\, \sum_{n=1}^{\infty}\, \frac{1}{n!} %
(\log (\alpha) \, s \wedge t )^n \quad s,t \in \mathbb{R},
\end{eqnarray*}
we see that even $ k_{\alpha} - 1 $ is positive definite, if $\alpha >1$.
\par
We denote further, for
$ \theta \, = \, e^{i \psi } \,  \mbox{ with } \psi \in [-\pi , \pi) , $
by
$ D_\theta $
the multiplication operator
\[ D_\theta \,:\, \mathcal{F} (\mathbb{R}_+) \, \rightarrow \, \mathcal{F} (\mathbb{R}_+) \,\]
defined by
\[ D_{\theta}f(t) \, = \, e^{i t \psi}f(t) \quad t \geq 0.\]
\par
As a theorem we state that the multiplication operators defined above act 
boundedly with respect to our sesquilinear form:
\begin{thm}
\label{minimum}
For all 
$ f \in \mathcal{F} (\mathbb{R}_+) $
we have
\[ \scalprod{D_{\theta}f}{D_{\theta}f} \, \leq \, C^2 \scalprod{f}{f}, \]
where 
$C \, = \, 1 + \, \frac{ 2\, \abs[ \psi ] }{ { \log(\alpha )}} .$
\end{thm}
{\it Proof: } 
Since for 
$ x \in \mathbb{R} $
\[ \alpha ^x \, = \, 1\,+\,\log{\alpha} \int_0^{x} \alpha ^u \, du ,\]
it follows that
\begin{eqnarray*} 
[ \, f \, , \, f \, ]_{\alpha} 
& \, = \, &\abs[\int_0 ^{\infty}\!f(t)\, dt] ^2 \,+\, \log{\alpha} \, \int_0 ^{\infty} \!\int_0 ^{\infty}\!
{\int_0 ^{s \wedge t}\!\alpha ^u \, du \, f(s)}\overline{f(t)}\, ds \, dt \\
               & \, = \, & \abs[\int_0 ^{\infty}\!\!f(t)\, dt] ^2 \,+\,\log{\alpha} \int_0 ^{\infty}\!\!\alpha ^u 
 \int_0 ^{\infty} \!\int_0 ^{\infty}\!\!{\chi_{[0,s \wedge t )}(u) f(s)}
\overline{f(t)} \, ds \, dt \, du \\
               & \, = \, & \abs[\int_0 ^{\infty}\!f(t)\, dt] ^2 \,+\,\log{\alpha} \int_0 ^{\infty}\!\alpha ^u
\abs[\int_u ^{\infty}\!f(t)\, dt] ^2 \, du.
\end{eqnarray*} 
Now, for arbitrary 
$ \eta > 0 $, $u\geq 0$:
\begin{eqnarray*} 
\abs[\int_u^{\infty} D_{\theta}f(t) \, dt] ^2 
& =  & \abs[\int_u^{\infty} \!e^{it\psi}f(t) \, dt] ^2 \\ 
& =  & \abs[\int_u^{\infty}\!\!
( e^{it\psi} -  e^{iu{\psi}})f(t) \, dt \, + \, \int_u^{\infty}\! \!e^{iu\psi}f(t) \, dt] ^2 \\
& \leq & (1 + \eta) \abs[\int_u^{\infty}\!\! f(t)\, dt]^2  + ( 1  +\frac{1}{\eta} )\abs[\int_u^{\infty}\!\!( e^{it\psi}-  e^{iu{\psi}})f(t) \, dt]^2.
\end{eqnarray*} 
The last integral can be estimated using an arbitrary
$ q \, > \, 1: $
\begin{eqnarray*} 
\abs[\int_u^{\infty}( e^{it\psi}\, - \, e^{iu{\psi}})f(t) \, dt]^2 
&  =  & \abs[\int_u^{\infty}i \psi \, 
\{ \int_u ^{t} e^{is\psi}\, ds \, \} \,f(t) \, dt\,]^2 \\
&  = & \abs[\psi] ^2 \abs[ \int_u ^{\infty} e^{is\psi} \, 
\{ \, \int_s ^{\infty} f(t) \, dt  \, \} \, ds]^2 \\
& \leq  & \abs[\psi] ^2 \int_u ^{\infty} q^{-2s} \, ds \, 
\int_u ^{\infty}q^{2s} \, \abs[\int_s ^{\infty}f(t)\, dt] ^2  \, ds \\
&  =  & \abs[\psi] ^2 \, \frac{1}{2 \log q} \, q^{-2u} \, 
\int_u ^{\infty}q^{2s} \, \abs[\int_s ^{\infty}f(t)\, dt]^2  \, ds. 
\end{eqnarray*} 
Hence, if 
$ \eta \, > \, 0 $
and
$ 1\,< \, q \, < \sqrt{\alpha}$
then
\begin{eqnarray*} 
\lefteqn{[  D_{\theta}f \, , \, D_{\theta}f  ]_{\alpha}  
\,  = \, \abs[ \int_0 ^{\infty} D_{\theta} f(t) \, dt]^2 \,+\,
\log {\alpha} \, \int_0^{\infty} \alpha ^u 
\abs[ \int_u ^{\infty} D_{\theta} f(t) \, dt]^2 \, du}&&\\
& \leq &
( 1 \, + \, \eta  )( \abs[ \int_0 ^{\infty} f(t) \, dt]^2 \,+\,\log {\alpha} \, \int_0^{\infty} \alpha ^u 
\abs[ \int_u ^{\infty}  f(t) \, dt]^2 \, du) \,+ \\
& & 
+ \, ( 1 \, + \, \frac{1}{\eta} )\abs[\psi]^2 \, 
\frac{1}{2 \log q}\,\left(\int_0 ^{\infty}q^{2s} \, \abs[\int_s ^{\infty}f(t)\, dt]^2
\, ds \,+\right.\\
& &\left. 
+  \, \log \alpha\int_0^{\infty}\! \alpha ^u q^{-2u}
{ \, \int_u^{\infty} \! q^{2s} \, 
\abs[\int_s ^{\infty}\!  f(t) \, dt]^2  \, ds \, } \,du\right) \\
& \, = \, &
( 1 \, + \, \eta ) \, [ \,f \,, \, f \, ]_{\alpha}\,
+ \, ( 1 \, + \, \frac{1}{\eta}) \, \abs[\psi]^2 \,
\frac{1}{2 \log q} \,\left(
\int_0 ^{\infty} \!  q^{2s}\,\abs[\int_s ^{\infty}\!  f(t) \, dt]^2 \,ds\,+\right.\\
& &\left.\int_0 ^{\infty} \!  q^{2s} \,\log \alpha \int_0^{s}\alpha ^u q^{-2u} \, du \,\abs[\int_s ^{\infty} \! f(t) \, dt]^2 \, \, ds \right)\\
&\, = \, & 
(1 \, + \, \eta ) \, [ \,f \,, \, f \, ]_{\alpha}\, + \\
& & 
+ \, ( 1 \, + \, \frac{1}{\eta} ) \, \abs[\psi]^2 \,
\frac{1}{2 \log q} \,
\int_0 ^{\infty} \!\! ( \, q^{2s}\,+\,\log \alpha 
\frac{\alpha^s  -  q^{2s}  }{ \log \alpha - 2 \log q} \, ) \, 
\abs[\int_s ^{\infty}\!\!  f(t) \, dt]^2  \,ds \\
& \, \leq \, &
(  1 \, + \, {\eta} )\, [ \,f \, ,\, f \, ]_{\alpha} \,+ \\
& & +\, ( 1 \, + \, \frac{1}{\eta} ) \, \abs[\psi]^2 \, \frac{1}{2 \log q} \,
\frac{\log \alpha}{ \log \alpha - 2 \log q} \, 
\int_0 ^{\infty}\! \alpha ^s 
\abs[\int_s ^{\infty} \! f(t) \, dt]^2 \,ds. 
\end{eqnarray*}  
Choosing
$ q \, = \, \alpha^{\frac{1}{4}} $
we have
\begin{eqnarray*} 
[ \,D_{\theta} f \, ,\,D_{\theta} f \, ]_{\alpha} & \, \leq \, & 
(\, 1 \, + \, \eta \, )  \, [\,f \,, \, f \, ]_{\alpha}\, +
\\
& & + \, (\, 1 \, + \, \frac{1}{\eta} \,) \, \abs[\psi]^2 \,
\frac{4}{ \log \alpha} \, \int_0 ^{\infty} \alpha ^s \, 
\abs[\int_s ^{\infty}  f(t) \, dt]^2 \,ds\\
& \, \leq \, & ( \, (\, 1 \, + \, \eta \, ) \, + \, (\, 1 \, + \, \frac{1}{\eta} \,)\, \abs[\psi]^2 \, \frac{4}{ \log^2 \alpha } \, ) \, [ \, f \, , \, f \, ]_{\alpha}.
\end{eqnarray*} 
By minimising this in 
$ \eta > 0 $
we obtain finally
\[ 
[ \, D_{\theta}f \, , \, D_{\theta}f \, ]_{\alpha}  \, \leq \, 
( \, 1 \, + \, \frac{2\,  \abs[\psi]}{ {\log \alpha }} \, )^2 \, [ \, f \, , \, f \, ]_{\alpha}.
\]
\par
We shall actually need to apply a discrete version of the above theorem, which
we want to formulate next.
To do this we define for finitely supported functions
\[g,f\,:\,\mathbb{Z}_+ \times \mathbb{Z}_+ \rightarrow \mathbb{C} \]
an, again positive definite, sesquilinear form  
\[ [\, f\,,\, g \,]_{\alpha} ^{\circ} \, = \scalproddis{f}{g} \]
and we let the multiplication operator 
$D_{\theta} $
be defined by 
\[ (D_{\theta}f)_n \, = \, \theta ^n f_n \quad n \in \mathbb{Z}_+ .\]
\par
\begin{thm}
\label{minimumdiscr}
For 
 all finitely supported functions
$ \quad f\,:\,\mathbb{Z}_+ \times \mathbb{Z}_+ \rightarrow \mathbb{C} :$
\[[\, D_{\theta}f\,,\, D_{\theta}f \,]_{\alpha} ^{\circ} \, \leq \,
\left( 1 + \, \frac{ 2\, \abs[ \arg \theta ] }{ { \log(\alpha )}}\right)
\, [\, f\,,\, f \,]_{\alpha} ^{\circ} .\]
\end{thm} 
{\it Proof:}
For a function $f:\mathbb{Z} \rightarrow \mathbb{C} $ and for $ n \in \mathbb{N}$ 
let  
$\Phi_n(f)$ be the function on $\mathbb{R}$ defined by
\begin{eqnarray*}
  \label{eq:embedding}
  \Phi_n(f)(x) & = & \sum_{k\in \mathbb{Z}}\, f_k\, n\, %
\chi_{[k,k+\frac{1}{n})}(x)%
\quad x\in \mathbb{R}.
\end{eqnarray*}
Here $\chi_{[a,b)}$ denotes the indicator function of the interval $[a,b)$.
\par
Then, if $f$ is supported on $\mathbb{Z}_+$:
\begin{eqnarray*}
  \label{eq:a}
  [\, f\,,\, f \,]_{\alpha} ^{\circ}& = & \lim_{n \rightarrow \infty}%
 [\, \Phi_n(f)\,,\, \Phi_n(f) \,]_{\alpha}.
\end{eqnarray*}
And, if $\theta = e^{i\psi}$, then
\begin{eqnarray*}
  \label{eq:b}
  \norm[1]{\Phi_n(D_{\theta}f)-D_{\theta}\Phi_n(f)} & \leq & %
\int_0^{\infty} \sum_{k\in \mathbb{Z}_+}\, \abs[f_k]\, \abs[e^{ik\psi}-e^{it\psi}]%
n\ \chi_{[k,k+\frac{1}{n})}(t)\, dt \\
& \leq & \sum_{k\in \mathbb{Z}_+}\, \abs[f_k]\, n\, \int_0^\frac{1}{n}%
\abs[1-e^{it\psi}]\, dt \\
& \leq & \sum_{k\in \mathbb{Z}_+}\, \abs[f_k]\, n\, \int_0^\frac{1}{n}%
\abs[\psi] \frac{1}{n} \, dt,
\end{eqnarray*}
which converges to zero as $n \rightarrow \infty$.
\par
Since for a finitely supported  $f$ all the above functions %
$\Phi_n(D_{\theta}f),\,  \Phi_n(f) \ldots$ have their support in a 
compact set on whose Cartesian product the kernel $k_{\alpha}$ remains 
bounded we obtain
\begin{eqnarray*}
  \label{eq:c}
\lim_{n \rightarrow \infty}  [\, \Phi_n(D_{\theta}f)\,,\, \Phi_n(D_{\theta}f) \,]_{\alpha}& = & \lim_{n \rightarrow \infty}%
 [\, D_{\theta}\Phi_n(f)\,,\, D_{\theta}\Phi_n(f) \,]_{\alpha}.
\end{eqnarray*}
Then, of course
\begin{eqnarray*}
  \label{eq:d}
  [\, D_{\theta}f\,,\, D_{\theta}f \,]_{\alpha} ^{\circ}& = & \lim_{n \rightarrow \infty}%
 [\, \Phi_n(D_{\theta}f)\,,\, \Phi_n(D_{\theta}f) \,]_{\alpha}\\
& = & \lim_{n \rightarrow \infty}%
 [\, D_{\theta}\Phi_n(f)\,,\, D_{\theta}\Phi_n(f) \,]_{\alpha}\\
& \leq & C_{\theta , \alpha} \, \lim_{n \rightarrow \infty}%
 [\, \Phi_n(f)\,,\, \Phi_n(f) \,]_{\alpha}\\
& = & C_{\theta , \alpha} \, [\, f\,,\, f \,]_{\alpha} ^{\circ}.
\end{eqnarray*}
\begin{rem}
Similarly to the above, one can extend Theorem~{\rm\ref{minimum}} to
finite Borel measures which are supported in  bounded subsets of $\mathbb{R}_+$.
\end{rem}
\begin{cor}
\label{cor:bdmult}%
If 
\[
\Lambda : \mathbb{Z}_+ \times \mathbb{Z}_+ \rightarrow \mathbb{C} 
\]
$(k,l) \mapsto \Lambda_{k,l}$ is a  positive definite kernel%
\footnotetext{positive definite only means positive semi-definite}
with finite support, then
with the same constant $C_{\theta,\alpha}$ as in the above theorems:
\[
\sum_{k,l} \heightheta^{k}
{\overline{\theta}^{l}}\alpha^{k \wedge l}\Lambda_{k,l}%
 \leq C_{\theta,\alpha}\sum_{k,l}\alpha^{k \wedge l}\Lambda_{k,l}.
\]
\end{cor}
{\it Proof:}
Since $\Lambda $ is positive definite with finite support, there exists 
a finite dimensional Hilbert space $\mathcal{H}$ and a finite sequence
$(\xi_k)_{k\geq 0}^{K}$ of its elements, such that
$ \Lambda_{k,l} = <\xi_k,\xi_l>$.
If $ (e_i)_{i\in I}$ is an orthonormal basis of $\mathcal{H}$, then\\[-1em]
\begin{eqnarray*}
\sum_{k,l}\heightheta^{k}%
{\overline{\theta}^{l}}\alpha^{k \wedge l}\Lambda_{k,l}
&=&
\sum_{i\in I}%
\sum_{k,l=0}^{K} \heightheta^{k}%
{\overline{\theta}^{l}}\alpha^{k \wedge l}<\xi_k,e_i><e_i,\xi_l>\\
&\leq&%
C_{\theta,\alpha}\sum_{i\in I}%
\sum_{k,l=0}^{K}\alpha^{k \wedge l}<\xi_k,e_i><e_i,\xi_l>\\
&=&C_{\theta,\alpha}\sum_{k,l}\alpha^{k \wedge l}\Lambda_{k,l}.
\end{eqnarray*}
In the rest of this section we state and prove, for later use, two
lemmas from Fourier analysis on the abelian groups $\mathbb{Z}$ and $\mathbb{Z}^k$.
\begin{lem}
\label{lem:gauss}
For $0<q<1$ there exists $\mu = \mu(q)>0$
such that 
\[
m \mapsto q^{\abs[m]^2 - \mu \abs[m]} \qquad m \in \mathbb{Z}
\]
is a positive definite function.
\end{lem}
{\it Proof:}
The orthogonality relations for the characters of 
the torus, the dual group of $\mathbb{Z}$ show that it suffices
to prove 
for sufficiently small $\mu>0$ the non-negativity of the Fourier
series
\[
\Psi_{\mu}(t)=%
\sum_{j\in \mathbb{Z}} q^{\abs[j]^2 - \mu \abs[j]} e^{-ijt} \qquad t \in
  [0,2\pi).
\]
To this end let $\tau = -\log q$ and notice that for $j \in \mathbb{Z}$:
\begin{eqnarray*}
q^{\abs[j]^2} &= & e^{-\tau j^2}\\
&=& \frac{1}{2\sqrt{ \pi \tau}} \int_{\mathbb{R}} e^{-\frac{x^2}{4\tau}}e^{ijx}\, dx\\
&=&\frac{1}{2\sqrt{ \pi \tau}} \int_{[0,2\pi)} %
\sum_{k \in \mathbb{Z}} e^{-\frac{(x+2 \pi k)^2}{4\tau}}e^{ijx}\, dx .
\end{eqnarray*}
Thus, by the uniqueness of Fourier inversion:
\[
\Psi_{0}(t)=
\sum_{j\in \mathbb{Z}} q^{\abs[j]^2}e^{-ijt} = %
\sqrt{\frac{\pi}{ \tau}}\sum_{k \in \mathbb{Z}} e^{-\frac{(t+2 \pi k)^2}{4\tau}} %
\qquad t \in [0,2\pi),
\]
which takes a strictly positive minimum $c_{\tau}>0$ on $[0,2 \pi)$. But
uniformly in $t \in[0,2 \pi)$:
\[
\abs[\Psi_{\mu}(t)-\Psi_{0}(t)] \leq \sum_{j\in \mathbb{Z}} %
\abs[q^{- \mu \abs[j]}-1] q^{\abs[j]^2},
\]
which for $\mu \searrow 0$  tends to $0$.
\par
The next lemma is prove similarly.
\begin{lem}
\label{lem:poisson}
For $0<q<1$ and $k \in \mathbb{N}$ there exists $\mu' = \mu'(q,k)>0$
such that 
\[
(m_1, \ldots ,m_k) \mapsto %
q^{( \abs[m_1]+\ldots+\abs[m_k]) - \mu'\abs[m_1+ \ldots + m_k]}, %
\qquad m =(m_1, \ldots ,m_k)\in \mathbb{Z}^k
\]
is a positive definite function on $\mathbb{Z}^k$.
\end{lem}
{\it Proof:}
This time we note first that
\[
\Phi(t)=\sum_{n=-\infty}^{\infty} q^{\abs[n]} e^{-int} = \frac{1-q^2}{1-2q\cos{t}+q^2}
\]
is strictly positive on the one dimensional torus $[0,2\pi)$.
\par
Hence on the $k$-dimensional torus
\[
\Psi_0(t_1, \ldots ,t_k)=\sum_{n_1=-\infty}^{\infty}\cdots \sum_{n_k=-\infty}^{\infty} %
q^{(\abs[n_1]+ \ldots + \abs[n_k])}e^{i(n_1t_1+ \ldots +n_kt_k)} %
= \prod_{l=1}^{k} \Phi(t_l)
\]
is strictly positive too. 
\par
Denote
 \[
\Psi_{\mu'}(t_1, \ldots ,t_k) = \sum_{n_1=-\infty}^{\infty}\cdots \sum_{n_k=-\infty}^{\infty} %
q^{(\abs[n_1]+ \ldots + \abs[n_k])-\mu'\abs[n_1+ \ldots +n_k]} %
e^{i(n_1t_1+ \ldots +n_kt_k)} .
\] 
Then, as $\mu' \searrow 0$
\[
\Psi_{\mu'}(t_1, \ldots ,t_k) \rightarrow \Psi_0(t_1, \ldots ,t_k)
\]
uniformly in $(t_1, \ldots ,t_k) \in [0,2\pi)^k$.
In fact, the absolute value of the difference is dominated uniformly by
\[
\sum_{n_1=-\infty}^{\infty}\cdots \sum_{n_k=-\infty}^{\infty}%
q^{(\abs[n_1]+ \ldots + \abs[n_k])}\left( q^{-\mu'\abs[n_1+ \ldots%
      +n_k]} - 1 \right).
\]
The only point to note is that this series is finite for small $\mu'>0$.
But, by the triangle inequality:
\[-\mu'\abs[n_1+ \ldots+n_k] \geq -\mu'(\abs[n_1]+ \ldots + \abs[n_k]),
\]
and since $q<1$:
\[
 q^{(\abs[n_1]+ \ldots + \abs[n_k])}q^{-\mu'\abs[n_1+ \ldots%
      +n_k]} \leq q^{(1-\mu)(\abs[n_1]+ \ldots + \abs[n_k])}.
\]
From this we infer for $\mu'<1$ the summability of the geometric series.
\section{Coxeter groups, their Cayley %
graph and their geometrical representation}
A group $G$ together with a finite generating set 
$S=\{s_1, \ldots ,s_n\}$
is called a Coxeter system, or Coxeter group,
if it has a presentation
\[
s_{1}^{2}, \ldots , s_{n}^{2},%
\quad (s_i s_j)^{m_{ij}} \qquad i,j \in \{1, \ldots ,n \} ,\, i \neq j.
\]
Here
$m_{ij}\in \{2,3, \ldots,\infty \}$
denotes the order of the product of the two generators $s_i,s_j \in S$,
$ i \neq j$,
setting $m_{ij}= \infty$ if
$s_i s_j $ has no finite order.
(If $s,s'$ are elements of $S$, $s=s_i$ and $s'=s_j$ for some 
$i,j \in \{1, \ldots n\}$ we shall write $\m{s}{s'}$ instead
of $m_{ij}$.)
The cardinality $n$ of $S$ is called the rank of the
Coxeter group. 
\par
The Cayley graph 
$\cayley{G}{S}$
is the graph with vertices just the group elements
$V=\{g: g \in G \}$.
Two vertices $u,v \in V$ are connected by an edge if
$u=vs$ for some $s\in S$.
We denote 
$E=\{ \{u,v\} : u=vs \mbox{ for some } s\in S\}$
the set of edges and notice that an edge 
connecting the vertices 
$g \in V \mbox{ and } gs  \in V$
is canonically labeled by a generator $s \in S $.
\par
We shall identify every edge with an image of the
unit interval 
$ [0,1] \subset \mathbb{R}$ 
obtaining thus a connected metric space $(\cayley{G}{S},{\rm \bf d})$.
A path in $\cayley{G}{S}$ then is a rectificable 
map $p:I \rightarrow \cayley{G}{S}$
defined on  some closed interval 
$I \subset \mathbb{R}$ 
into $ \cayley{G}{S}$.
\par
If we denote by 
${\rm \bf l}:G \rightarrow \mathbb{N}$
the length function with respect to the generating set, i.e.\ if:
\[
\la{g} = \inf \{m: g=s_{i_1}\cdot \ldots \cdot s_{i_m},%
 \, s_{i_j} \in S \cup S^{-1} \} \qquad g \in G
\]
then clearly this relates to the distance in $\cayley{G}{S}$ by:
\[
\di{g}{h} = \la{g^{-1}h} \qquad g,h \in G.
\]
It is obvious from the above definitions that
the action of the group $G$ on itself by left multiplication,
$g: h \mapsto gh, \quad g,h \in G$,
extends to an action of $G$ by isometries of the metric space $(\cayley{G}{S},d)$.
\par
In an equation $ g=w_1 \ldots w_m$, we shall call the right hand side
a reduced representation of $g$, if $w_i \in S , \; i =1, \ldots, m$
and $\la{g}=m$. A product $u_1 \ldots u_k$ will be called  reduced
if $\la{u_1 \ldots u_k} = \la{u_1}+ \ldots +\la{u_k}$.
The void word represents the identity $\e$ of $G$.
\par 
A useful tool for a finitely generated Coxeter group is its 
representation as a discrete subgroup of the general linear
group of a finite dimensional real vector space $E$ of dimension $\card{S}$
(see chapter V of \cite{Bourbaki_68}). We shall denote it:
\[
\sigma : G \rightarrow {\rm Gl}(E).
\]
A self inverse element $t \in G, \; t \neq \e$ will be called a reflection
if $\sigma(t)$, or equivalently $\sigma^{\ast}(t)$,
is a reflection. The corollaire in 3.2, chapter V of \cite{Bourbaki_68},
asserts that any reflection is conjugate to some generator.
($T=\{g^{-1}sg \ : \ s\in S,\ g\in G \}$ will denote the set of all 
reflections). 
We shall denote
for $g\in G$:
\begin{eqnarray}
N_g & = &\{ t\in T: \la{tg}<\la{g}\}
\end{eqnarray}
\par
\begin{rem}
\begin{description}
\item[~~(i)]
\label{thm:length}%
The length of a group element
$g \in G$ is given by:
\[ 
\la{g} = \card \{ t \in T : \la{tg}<\la{g}\}.
\]
\item[~(ii)]
If $ g = w_1 \ldots w_n$ is a reduced representation, 
then
\begin{eqnarray}
\label{eq:roots}%
N_g & = & \{ w_1,\, w_1w_2w_1, \, 
\ldots , \, w_1\ldots w_{n-1} w_n
w_{n-1}\ldots  w_1 \} .
\end{eqnarray}
\end{description}
\end{rem}
\par
The next  theorem is a reformulation, in our setting, of one
of~\cite{BozJanSpatz_88}.
\begin{thm}
\label{thm:boz}%
For $g,h \in G$:
\[
\di{g}{h} = \card N_g \bigtriangleup N_h = %
\sum_{t\in T} \abs[\chi_{N_g}(t) - \chi_{N_h}(t)]^2.
\]
\par
Hence $g \mapsto \la{g}$ is a negative definite function.
\end{thm}
We shall actually need some consequences
of the action of $G$ by means of the contragradient representation
\[
\sigma^{\ast} : G \rightarrow {\rm Gl}(E^{\ast}).
\]
on the Tits cone $U$, see $\mbox{n}^{\mbox{\small o}}$ 4.6. of \cite{Bourbaki_68}
\begin{rem}
\begin{description}
\item[~~(i)]
If for $t,t' \in T$ the hyperplanes stabilised by $\sigma^{\ast}(t)$,
respectively $\sigma^{\ast}(t')$, intersect inside the Tits cone $U$, then 
the product $tt'$ has finite order. 
(This is contained in exercise 2c) and 2d) of {\rm \cite{Bourbaki_68}})
\item[~~(ii)]
$G$ contains a normal torsion free subgroup $\Gamma$ of finite index.
(see  loc.\ cit.\ exercise 9).)
\end{description}
\end{rem}
\par
Now $G$, and hence $\Gamma$, act on the set of reflections
$T$ by conjugation. Let
\[
T=T_1 \dot{\cup} T_2 \dot{\cup} \ldots \dot{\cup} T_{\Lambda}
\]
be a decomposition into $\Gamma$--orbits. 
\begin{lem}
\label{lem:disjoint}
If $t,t'\in T$ belong to the same $\Gamma$--orbit,
then either $t=t'$ or the hyperplanes stabilised
by $\sigma^{\ast}(t)$, and
 $\sigma^{\ast}(t')$ respectively, do not intersect inside the Tits cone $U$
\end{lem}
\begin{proof}
If we assume that the stabilised hyperplanes
intersect inside $U$, then $tt'$ has a finite order.
On the other hand, if $t$ and $t'$ both belong to the same $\Gamma$--orbit,
then, for some $t_0\in T$ and some $\gamma,\gamma'\in \Gamma$:
\[
t=\gamma^{-1}t_0\gamma \mbox{ and } t'={\gamma'}^{-1} t_0 {\gamma'}
\]
Since $t_0^2=\e$, and $\Gamma$ is normal
\[tt'= \gamma^{-1}t_0\gamma{\gamma'}^{-1} t_0 {\gamma'} \in \Gamma.
\]
Because $\Gamma$ is torsion--free it follows that this
element of finite order equals the identity.
\end{proof}
\par
For $g\in G$ we momentarily fix a reduced decomposition 
and order $N_g$ according to (\ref{eq:roots}). 
Further we endow the sets
\[
N_{g}^i=N_g \cap T_i, \quad i=1, \ldots, \Lambda
\]
with the order inherited as  subsets of $N(g)$. In fact, the order
obtained on the subsets $N_{g}^i$ does not depend on the
reduced decomposition chosen.
\begin{lem}
Let $g\in G$ be given. Then for any $u\in G$ and for $i\in\{1,\ldots,\Lambda\}$
\[N_{g}^i \cap N_u\]
is an initial segment of $N_{g}^{i}$
\end{lem}
\begin{proof}
We have to show, that 
$t\in N_{g}^i \cap N_u$
and 
$t'\in N_{g}^i,\ t'<t$
imply 
$t'\in N_u$.
\par 
Denote $H_t$, $H_{t'}$ the hyperplanes in $E^*$
fixed by $\sigma^{\ast}(t)$, respectively by
$\sigma^{\ast}(t')$.
Then $t\in N_g$ means, that $H_t$ separates the
fundamental chamber $C$ from $\sigma^{\ast}(g)C$,
similarly for $u$.
Moreover $t'<t$ shows, that there is a point in
$H_t \cap U$ separated from $C$ by $H_{t'}$.
Since $H_t$ and $H_{t'}$ do not intersect inside $U$, we
conclude that all of $H_t\cap U$ and $C$ ley on different sides
of $H_{t'}$.
Since $\sigma^{\ast}(g)C$ and $\sigma^{\ast}(u)C$
are separated from $C$ by $H_t$ we conclude that any
line segment from $\sigma^{\ast}(u)C$ to $C$ must
intersect $H_{t'}$, meaning that $t'\in N_u$.
\end{proof}
\par
Let us denote 
\begin{eqnarray*}
\label{eq:pos}
N^g = \{ N_g \cap N_u : \, u \in G \}.
\end{eqnarray*}
\begin{prop}
\label{prop:posdef}
For $0<r<1$ there exists a constant $\kappa$ such that 
\[
(U,V) \mapsto r^{\kappa (\card U \wedge \card V) %
 + \card ( U \bigtriangleup V)}
\]
is a positive definite kernel on $N^g\times N^g$.
\end{prop}
{\it Proof:}
For $U=N_u\cap N_g \in N^g$ we denote
$U_j=T_j\cap U$, $j=1,\ldots, \Lambda$;
 similarly for $V \in N^g$.
Further denote $u_j=\card U_j$ and $v_j=\card V_j$,
the respective cardinalities.
\par
Lemma~\ref{lem:poisson} provides a constant $\mu'>0$, such that
\[
(U,V) \mapsto r^{-\mu' \abs[\sum u_j-v_j] + \sum \abs[u_j-v_j]}
\]
is positive definite. Since
\begin{eqnarray*}
2(\card U \wedge \card V) &= & \card U +   \card V - \abs[\card U -  \card V]\\
&=&  \card U +   \card V - \abs[ \sum (u_j-v_j)],
\end{eqnarray*}
and since 
\[
(U,V) \mapsto r^{\mu'(\card U +   \card V)} \; = \; %
r^{\mu'\card U} \cdot r^{\mu'\card V} 
\]
is positive definite,
it suffices to show
\[
 \card (U \bigtriangleup V) = \sum \abs[u_j-v_j] .
\]
Because then 
\begin{eqnarray*}
r^{2{\mu'} (\card U \wedge \card V) + \card ( U \bigtriangleup V)} %
&=&
r^{\mu'(\card U +   \card V)} \cdot %
r^{-\mu'  \abs[\sum u_j-v_j] + \sum\abs[u_j-v_j]} 
\end{eqnarray*}
is positive definite as a product of positive definite kernels.
Now,
\begin{eqnarray*}
\card (U \bigtriangleup V)&=&\sum_{t\in N_g} \abs[\chi_U (t)- \chi_V(t)]^2\\
&=&\sum_j \sum_{t\in N_g\cap T_j} \abs[\chi_U (t)- \chi_V(t)]^2.
\end{eqnarray*}
Because $\chi_U - \chi_V$ takes only values in $\{-1,0,1\}$ we can omit the
squares. Moreover, since the sets $U_j$ and $V_j$ 
contain all predecessors of their elements,
$\chi_{U_{j}}-\chi_{V_{j}}$ is either a non-negative
or a non-positive function. Thus we may continue:
\begin{eqnarray*}
\card (U \bigtriangleup V)&=&\sum_j\sum_{t\in N_g} %
\abs[\chi_{U_j} (t)- \chi_{V_j}(t)]\\
&=& \sum_j \abs[\sum_{t\in N_g}\chi_{U_j} (t)- \chi_{V_j}(t)]\\
&=& \sum_j \abs[u_j-v_j].
\end{eqnarray*}
\section{Construction of a series of representations}
Theorem~\ref{thm:boz} shows that for any Coxeter group
\[
\di{.}{.} : G \times G \rightarrow \mathbb{Z}_+
\]
is a negative definite kernel on $G$. By Schoenberg's
theorem (see e.g.\ \cite{BergForst_75}), then,
for $r \in (0,1)$,
\[ 
g \mapsto r^{\la{g}} \qquad g \in G
\]
is a positive definite function.
\par 
We are going to consider modifications of the Gelfand--Naimark--Segal
representation constructed from this positive definite functions. So let 
\begin{eqnarray*}
\mathcal{F} & = & \{ f:G \rightarrow \mathbb{C} \; \mbox{ with finite support } \}\\
<f,h>_{r} & = & \sum_{u,v \in G} r^{\di{u}{v}} f(u) \overline{h(v)}%
 \qquad f,h \in \mathcal{F}\\
\norm[r]{f} & = & (<f,f>_{r})^{\frac{1}{2}} \qquad f \in \mathcal{F}\\
\pi_r (g) f(u) & = & f(g^{-1}u) \qquad g, u \in G, \, f \in \mathcal{F}.
\end{eqnarray*}
\begin{prop}
The kernel of $ \norm[r]{.} \mbox{ on } \mathcal{F}$ equals $\{0\}$.
\end{prop}
{\it Proof:}
We shall show, that for a finite set $u_1,\ldots,u_n \in G$ 
the matrix 
\[
\left(r^{\di{u_i}{u_j}}\right)_{i,j=1}^n
\]
is non-degenerate.
Denote $\tau = - \log r$ and consider 
the functions $\chi_i =\sqrt{\tau}\chi_{N_{u_i}}$ as
elements of the Hilbert space $l^2(T)$. 
By Theorem~\ref{thm:boz} 
\begin{eqnarray}
r^{\di{u_i}{u_j}}&=&\exp(\norm{\chi_i - \chi_j}^2)\\
& =&\label{eq:nondeg}%
 \exp(\norm{\chi_i}^2) \exp(-2<\chi_i,\chi_j>)\exp(\norm{\chi_j}^2).
\end{eqnarray}
Proposition~2.2 of~\cite{Guich_72} implies that
\[
 \left( \exp(-2<\chi_i,\chi_j>)\right)_{i,j=1}^n
\]
is non-degenerate.
But then, from (\ref{eq:nondeg}), we see that 
$ \left(r^{\di{u_i}{u_j}}\right)_{i,j=1}^n$
is non-degenerate too.
\par
For a given $\theta \in \mathbb{C},\, \abs[\theta]=1$ of absolute value one we
define an equivariant cocycle
\[
c_{\theta} : \, G \times G \rightarrow \mbox{\rm Gl}(\mathcal{F})
\]
for $\pi_r$ by:
\[
c_{\theta}(u,g)f\, (v) = \theta^{(\la{u^{-1}v} - \la{g^{-1}v})}f(v) \qquad f
\in \mathcal{F}, \; u,v,g \in G.
\]
For any $r \in (0,1)$ the cocycle equalities
\begin{eqnarray*}
c_{\theta}(u,u) & = & {\rm id},\\
c_{\theta}(u,g)c_{\theta}(g,v)&=&c_{\theta}(u,v),
\end{eqnarray*}
and the equivariance
\begin{eqnarray*}
\pi_r (g)\, c_{\theta}(v,u) & = & c_{\theta} (gv,gu) \, \pi_r (g)
\end{eqnarray*}
are easily checked.
\par
Thus for $z \in \mathbb{C}$ with $\abs[z] < 1$, which we write
as $z=\theta r,\, r\in (0,1),\ \abs[\theta]=1$, we may define a representation 
of $G$ on $\mathcal{F}$ by
\begin{eqnarray*}
\pi_z (g) f & = &  c_{\theta}(\e,g)\circ \pi_r (g) f \qquad f \in \mathcal{F}, \; g \in G.
\end{eqnarray*}  
\par
We are about to formulate a criterion for the boundedness of this 
representation on the semi--normed space $\mathcal{F}, \, \norm[r]{.}$. To do this 
let me denote, for $g \in G$:
\begin{eqnarray*}
\label{eq:posdefi}
N^g = \{ N_g \cap N_u : \, u \in G \}.
\end{eqnarray*}
This is a subset of the potency set of $N_g$. In the cases which are of
interest to us its cardinality will be much smaller than $2^{\la{g}}$.
\par
\begin{thm}
\label{thm:ub}
Let $(G,S)$ be a Coxeter group.
Then for
$ z=\theta r \in D = \{ z \in \mathbb{C} : \abs[z] <1\}$ there exists a 
uniformly
bounded representation $ ( \pi_z, H_r) $ such that for some $\xi_0 \in H_r$:
\[
<\pi_z(g) \xi_0, \xi_0>_r = z^{\la{g}} \qquad g \in G.
\]
Moreover, for some constant $\kappa$, depending only on $(G,S)$: 
\[
\sup_{g\in G}\norm{\pi_z(g)} \leq %
1 + \frac{2 \abs[\arg(z^2)] }{\kappa { \abs[\log r]}}.
\] 
\end{thm}
{\it Proof:} Since $\pi_r (g) : G \rightarrow \mbox{\rm Gl}(\mathcal{F})$ is
a representation by invertible isometries, it suffices to show, that
$c_{\theta}(e,g)$ is bounded by 
$ 1 + \frac{2 \abs[\arg(z^2)] }{\kappa { \abs[\log r]}}$.
\par
We note that for $u,v \in G $
\begin{eqnarray*}
\di{u}{v} & = &  \sum_{t\in N_g} \abs[ \chi _{N_u}(t) - \chi _{N_v}(t)]^2
+ \sum_{t\in N_{g}^c} \abs[ \chi _{N_u}(t) - \chi _{N_v}(t)]^2 \\
& = & \card ((N_g \cap N_u) \bigtriangleup (N_g \cap N_v)) + d'(u,v),
\end{eqnarray*}
where 
$N_{g}^c$ denotes the complement of $N_{g}$ in the set of all reflections and
\begin{eqnarray*}
d'(u,v) & = &  \sum_{t\in N_{g}^c} \abs[ \chi _{N_u}(t) - \chi _{N_v}(t)]^2 
\end{eqnarray*}
is a negative definite kernel on $G$.
Thus for  $f \in \mathcal{F}$
\begin{eqnarray*}
(U,V) &\mapsto & \sum_{u \in U, v \in V} r^{d'(u,v)} f(u) \overline{f(v)}
\end{eqnarray*}
is positive definite on $N^g$.
For $ k \in \mathbb{Z}_+$ let $ E_k = \{ U \in N^g : \, \card U = k\}$. 
From the assumption on $g$ and $\kappa$ it follows, using
Schur's theorem (see e.g.\ \cite{BergForst_75}), that
\begin{eqnarray*}
(k,l) \, \mapsto \, %
\Lambda_{(k,l)} & = & \sum_{U\in E_k, V\in E_l}%
 r^{\kappa \card U \wedge \card V + \card ( U \bigtriangleup V)}%
\sum_{u \in U, v \in V} r^{d'(u,v)} f(u) \overline{f(v)}
\end{eqnarray*}
is a positive definite kernel on $\mathbb{Z}_+$.
In Theorem~\ref{minimumdiscr} we let $\alpha = r^{-\kappa}$ and conclude
\begin{eqnarray*}
\sum_{k,l=0}^{\infty} \heightheta^{2k}
{\overline{\theta}^{2l}} {\alpha^{k \wedge l}} \Lambda_{(k,l)} & %
\leq & C_{\theta^2,\alpha}%
\sum_{k,l=0}^{\infty}%
{\alpha^{k \wedge l}} \Lambda_{(k,l)}.
\end{eqnarray*}  
The  left hand side of this inequality computes to 
$\norm[r]{c_{\theta}(\e,g)f}^2$ and the right
one to $C_{\theta, r^{-\kappa}}\norm[r]{f}^2$. 
\par
In fact, we notice that 
\[
\la{u}+\la{g}-\la{g^{-1}u}=%
\card  N_u+\card N_g -\card N_g \bigtriangleup N_u =%
 2\card (N_g \cap N_u)
\]
and compute:
\begin{eqnarray*}
\lefteqn{\sum_{k,l=0}^{\infty}{\heightheta^{2k}}{\overline{\theta}^l} %
{\alpha^{k \wedge l}}\Lambda_{(k,l)} \; = }& & \\
&=&\sum_{k,l=0}^{\infty}{\heightheta^{2k}}{\overline{\theta}^{2l}} %
r^{-\kappa({k \wedge l})}\sum_{U\in E_k, V\in E_l}%
 r^{\kappa (\card U \wedge \card V)} r^{ \card ( U \bigtriangleup V)}%
\sum_{u \in U, v \in V} r^{d'(u,v)} f(u) \overline{f(v)}\\
&=&\sum_{k,l=0}^{\infty}{\heightheta^{2k}}{\overline{\theta}^{2l}} %
\sum_{U\in E_k, V\in E_l}\sum_{u \in U, v \in V}%
r^{ \card ((N_g \cap N_u) \bigtriangleup (N_g \cap N_v)) +d'(u,v) } %
f(u) \overline{f(v)}\\
&=&\sum_{k,l=0}^{\infty}{\heightheta^{2k}}{\overline{\theta}^{2l}} %
\sum_{u\in E_{k}', v\in E_{l}'} r^{\di{u}{v}} f(u) \overline{f(v)}\\
&=&\sum_{u\in G, v\in G} \heightheta^{2\card ((N_g \cap N_u))} %
{\overline{\theta}^{2\card ((N_g \cap N_v))}} %
r^{\di{u}{v}}f(u) \overline{f(v)}\\
&=&\sum_{u\in G, v\in G}\heightheta^{(\la{u}-\la{g^{-1}u})} %
{\overline{\theta}^{(\la{v}-\la{g^{-1}v})}} %
r^{\di{u}{v}}f(u) \overline{f(v)},
\end{eqnarray*}
where, for $k \in \mathbb{Z}_+$, we denoted %
$E_{k}' = \{ u \in G : \card (N_g \cap N_u)=k\}$.
The right hand side is computed similarly and
checking for the
constant finishes the proof.
\par
We shall denote $H_r$ the Hilbert space
obtained by completing $\mathcal{F}$ with respect to $ \norm[r]{.}$.
The above gives the announced bound on the norm of the operators 
$\pi_z(g),\, z\in D,\, g \in G $.
\par
If we denote for $u\in G$ by $\delta_u$ the point mass one at $u$, then
\begin{eqnarray*}
\pi_z(g)\delta_u\; (x)
&=& c_{\theta}(\e,g)%
\delta_u(g^{-1}.)\;(x)\\
&=&\theta^{\la{x}-\la{g^{-1}x}}\delta_u(g^{-1}x)\\
&=&\theta^{\la{gu}-\la{u}}\delta_{gu}(x) \qquad x\in G.
\end{eqnarray*}
Hence,
\begin{eqnarray*}
<\pi_z(g)\delta_{\e},\delta_{\e}>_r%
&=&\theta^{\la{g}}<\delta_{g},\delta_{\e}>_r\\
&=&\theta^{\la{g}}r^{\la{g}}.
\end{eqnarray*}
\begin{rem}\rm
\label{rem:ub}%
Assume that the conditions of the corollary are satisfied.
\begin{description} 
\item[(1)]
For any $u,v\in G$ the map $\theta \mapsto c_{\theta}(u,v)$ 
is a group homomorphism from the circle group (the torus) into
the bounded invertible operators on $H_r$.
It is continuous for the strong operator topology on $B(H_r)$.
\item[(2)]
If $\chi: g \mapsto (-1)^{\la{g}}$ , then $\chi$ is a character of the
Coxeter group $G$. For $z\in D$ the tensor product representation
$\chi \otimes \pi_z$ is canonically isomorphic to $\pi_{-z}$.
\item[(3)]
Complex conjugation of functions in $\mathcal{F}$ defines a conjugate linear
intertwining operator between $\pi_z$ and $\pi_{\overline{z}}$.
\end{description}
\end{rem}
\section{Weak amenability}
To prove the weak amenability of $G$ we can not directly apply the
Theorem of \cite{Valette_93}, since our
series of representations is not realized on one Hilbert space.
But in the proof of the mentioned theorem it is only used
that $z \mapsto \varphi_z$ is analytic as
a function from $D$ to the completely bounded
multipliers of $A(G)$.
This can be proved by a method of Pytlik and Szwarc \cite{PytSzw_86},
see also \cite{Szwarc_91}.
\par
For this
we denote by $\chi_n$  the characteristic function 
of the set $\{ g\in G \;: \, \la{g}=n\}$ of group elements of length $n\in \mathbb{N}$.
\begin{prop}
\label{prop:char}%
For some constant $\kappa$, depending only on $(G,S)$: 
\[
\sup_{g\in G}\norm[M_0A(G)]{\chi_n} \leq %
2\pi\; e\; (1 + \frac{4\pi}{ \kappa} n).
\]
\end{prop}
{\it Proof:}
\par
Let for $0<r<1$ denote $\tilde{\pi}_{r}:G \rightarrow B(L^2([0,2\pi],H_r))$
the direct integral representation on $L^2([0,2\pi],H_r)$ defined by:
\[
(\tilde{\pi}_{r}(g)f)\, (t) = \pi_z(g)(f(t)),\; \mbox{where }z=re^{it},\;%
f\in L^2([0,2\pi],H_r),\;g\in G.
\] 
Now, let $f,h\in L^2([0,2\pi],H_r)$ denote the square integrable $H_r$
valued functions: 
\[f: t \mapsto \delta_{\e} \mbox{ and } h: t \mapsto e^{int}\delta_{\e}.
\]
Then, for any $g\in G$:
\begin{eqnarray*}
<\tilde{\pi}_{r}(g)f,h> &=& \int_{0}^{2\pi}%
<\pi_{re^{it}}\delta_{\e},\delta_{\e}>_r%
e^{-int} \;dt\\
&=&\int_{0}^{2\pi} r^{\la{g}} e^{i\la{g}t}e^{-int} \;dt\\
&=&%
\left\{%
\begin{array}{lcl}
2\pi\; r^{\la{g}}& \mbox{ if }& n= \la{g},\\
0& \mbox{ if } & n \neq \la{g}
\end{array}\right.\\
&=& 2\pi\; r^n \chi_n(g).
\end{eqnarray*}
Since $\norm{f}=\norm{h}=\sqrt{2\pi}$, we infer from Theorem~\ref{thm:ub}:
\begin{eqnarray*}
\norm[M_0A(G)]{\chi_n} &\leq& r^{-n}\sup_{g\in G} \sup_{\abs[\theta]=1}%
\norm{\tilde{\pi}_{r\theta}(g)}\; \norm{f}\; \norm{h}\\
&\leq&2\pi\; r^{-n}(1 + \frac{4\pi}{ \kappa { \abs[\log r]}}).
\end{eqnarray*}
Here we may take $r=e^{-\frac{1}{n}}$ on the right hand side and obtain the
constant given in the statement.
\begin{cor}
The function $z \mapsto z^{\la{.}}$ is analytic on $D$. 
\end{cor}
{\it Proof:}
We just note that the series
\[
z^{\la{.}} = \sum_{n=0}^{\infty} \chi_n \, z^n
\]
is norm convergent in $M_0A(G)$, on the whole open disk $D$.
\par
Arguing either as in \cite{Valette_93}\ or as in the proof
of Theorem 6 in \cite{Szwarc_91}, see also the article
of de Canni\`{e}re and Haagerup \cite{CanHaa_84}, we obtain:  
\begin{thm}
A Coxeter group $(G,S)$ is weakly amenable with
Cowling--\-Haa\-gerup constant one.
\end{thm}
The proof of the theorem is immediate from the following lemma, which we state
and prove for the readers convenience.
\begin{lem}
\label{lem:weaklyamenable}%
Let $G$ be a locally compact group, and $z \mapsto \varphi_z$
an analytical  map from the unit disk $D=\{z\in \mathbb{C} : \abs[z]<1\}$
to $M_0A(G)$, such that
\begin{description}
\item[~~(i)]
$\displaystyle \varphi_z = \sum_{n=0}^{\infty} \psi_n\ z^n $,
fore some $\psi_n \in A(G)$,
\item[~(ii)]
$\varphi_r$ is an element of the 
unit sphere of the Fourier--Stieltjes algebra, for $r\in [0,1)$, 
\item[(iii)]
locally uniformly on $G$:
\[\lim_{r\rightarrow 1} \varphi_r = 1.\]
\end{description}
Then $G$ is weakly amenable with Cowling--Haagerup constant $1$.
\end{lem}
\begin{proof}
From theorem $\rm B_2$ of \cite{GranLei_79} we infer that
for all $\psi\in A(G)$:
\begin{eqnarray}
\label{eq:multconv}%
\norm[A(G)]{\varphi_r \psi -\psi} &\rightarrow 0 ,& \mbox{ as } r\rightarrow 1.
\end{eqnarray}
Now taking the Fejer kernel on the torus $\mathbb{T}$:
\[
F_N(e^{it}) = \sum_{\abs[k]\le N } \left( 1- \frac{\abs[k]}{N+1}\right) e^{ikt},
\]
we have
\begin{description}
\item[1)]
$F_N \geq 0 $ and $\frac{1}{2\pi} \int_0^{2\pi} F_N(e^{it}\ dt) =1$,
\item[2)]
for all $f\in C(\mathbb{T})$:
$F_N \ast f (\theta) \rightarrow f (\theta),\quad \forall \theta \in \mathbb{T}$ as
$N \rightarrow \infty$.
\end{description}
Let
\[ \psi_{N,r} = \frac{1}{2\pi} \int_0^{2\pi} F_N(e^{it}) \varphi_{re^{it}} dt.\]
Then $\psi_{N,r} \in M_0A(G)$ and 
\begin{eqnarray}
\label{eq:fejer}%
\norm[M_0A(G)]{\psi_{N,r}-\varphi_{r}}&\leq&
\frac{1}{2\pi} \int_0^{2\pi} F_N(e^{it})
\norm[M_0A(G)]{\psi_{N,r}-\varphi_{re^{it}}} dt,
\end{eqnarray}
which, since $z \mapsto \varphi_z$ is continuous on $D$, converges
to $0$ as $N\rightarrow \infty$, by {\bf 2)}.
\par
On the other hand $\psi_{N,r}\in A(G)$, since
\begin{eqnarray*}
\psi_{N,r}(g)&=&%
\frac{1}{2\pi} \int_0^{2\pi} F_N(e^{it}) \sum_n\psi_n(g)\ {(re^{it})}^n dt\\
&=&\sum_{n=0}^{N}(1-\frac{n}{N+1})\psi_n(g).
\end{eqnarray*}
Now from (\ref{eq:multconv}) and (\ref{eq:fejer}) it is easy to construct
an approximate unit in $A(G)$ with its completely bounded multiplier norm
bounded by $1$.
\end{proof}

%
%
\vspace{\fill}
{{\bf Author's address}\\
Gero Fendler\\
Finstertal 16\\
D-69514 Laudenbach\\
Germany\\
{\small e-mail: fendler@math.uni-sb.de}}
\end{document}